\documentclass[11pt]{article}

\usepackage{relsize}
\usepackage{phonetic}
\usepackage{boldline}
\usepackage{bm}
\usepackage{epsfig}
\usepackage{caption}
\usepackage{amsmath,amsfonts} % Typical maths resource packages
\usepackage{graphics}                 % Packages to allow inclusion of graphics
\usepackage{graphicx}
\usepackage{empheq}
\usepackage{array}
\usepackage{epsf}
\usepackage{color}                % For creating coloured text and background
\usepackage{accents}
\usepackage{tcolorbox}
\usepackage{appendix}
\usepackage{float}
\usepackage{hyperref}
\usepackage{multirow}
\usepackage{array}
\usepackage{booktabs}
\usepackage{mathrsfs}

\newtheorem{theorem}{Theorem}[section]
\newtheorem{lemma}[theorem]{Lemma}

\newtheorem{corollary}[theorem]{Corollary}
\newtheorem{example}{Example}[section]

\newenvironment{proof}[1][Proof]{\begin{trivlist}
		\item[\hskip \labelsep {\bfseries #1}]}{\end{trivlist}}

\newcommand{\nn}{\nonumber}

\def\refe#1{(\ref{#1})}
\def \bb#1{\setbox0=\hbox{$#1$}
	\kern-.025em\copy0\kern-\wd0 \kern.05em\copy0\kern-\wd0
	\kern-.025em\raise.0433em\box0}
\def \endproof{\vrule height8pt width 5pt depth 0pt}

\newcommand{\qed}{\nobreak \ifvmode \relax \else
	\ifdim\lastskip<1.5em \hskip-\lastskip
	\hskip1.5em plus0em minus0.5em \fi \nobreak
	\vrule height0.75em width0.5em depth0.25em\fi}

\newcommand{\bs}{{\bm{\sigma}}}

\newcommand{\bv}{{\bm{v}}}

\newcommand{\bH}{{\mathbf{H}}}
\newcommand{\bn}{{\mathbf{n}}}
\newcommand{\D}{\mathrm{div}}

\newcommand{\n}{\mathbf{n}}

\newcommand{\ee}{\epsilon}

\newcommand{\RT}{\mathcal{R\! T}}
\newcommand{\DG}{\mathcal{D G}}
\newcommand{\x}{\bm{x}}
\renewcommand{\P}{\mathcal{P}}
\begin{document}
	
\title{Analysis of A Mixed Finite Element Method for Poisson's Equation with Rough Boundary Data}

\author{
	\setcounter{footnote}{0}
	Huadong~Gao,
	\footnote{
		School of Mathematics and Statistics,
		Huazhong University of Science and Technology,
		Wuhan 430074, P.R. China
		({\tt huadong@hust.edu.cn}).
		The work of the author was supported in part by the
		National Science Foundation of China No. 12231003
		and National Key Research and Development Program of China (No. 2023YFC3804500).}
	\quad
	Yuhui~Huang
	\footnote{
		School of Mathematics and Statistics,
		Huazhong University of Science and Technology,
		Wuhan 430074, P.R. China.
		({\tt m202370045@hust.edu.cn})
	}
	~~ and ~~
	Wen~Xie
	\footnote{
		School of Mathematics and Statistics,
		Huazhong University of Science and Technology,
		Wuhan 430074, P.R. China.
		({\tt w\_xie@hust.edu.cn})
	}
}
\date{}

\maketitle

\begin{abstract}
  This paper is concerned with finite element methods for Poisson's equation with rough boundary data. Conventional methods require that the boundary data $g$ of the problem belongs to $H^{1/2} (\partial \Omega)$. However, in many applications one has to consider the case when $g$ is in $L^2(\partial \Omega)$ only. To this end, very weak solutions are considered to establish the well-posedness of the problem. Most previously proposed numerical methods use regularizations of the boundary data. The main purpose of this paper is to use the Raviart--Thomas mixed finite element method to solve the Poisson equation with rough boundary data directly. We prove that the solution to the proposed mixed method converges to the very weak solution. In particular, we prove that the convergence rate of the numerical solution is $O(h^{1/2})$ in convex domains and $O(h^{s-1/2})$ in nonconvex domains, where $s > 1/2$ depends on the geometry of the domain. The analysis is based on a regularized approach and a rigorous estimate for the corresponding dual problem. Numerical experiments confirm the theoretically predicted convergence rates for the proposed mixed method for Poisson's equation with rough boundary data.
	

\end{abstract}

\section{Introduction}
\setcounter{equation}{0}
In this paper,
we consider the Poisson equation with Dirichlet boundary condition
\begin{empheq}[left=\empheqlbrace]{align}
	& - \Delta u = f
	&& \text{in} \quad   \Omega,
	\label{pde}
	\\
	& u = g
	&& \text{on} \quad \Gamma \! ,
	\label{pde-bc}
\end{empheq}
where $\Omega$ is a bounded Lipschitz polygonal/polyhedral domain in $\mathbb{R}^d(d=2,3)$
and $\Gamma$ denotes the boundary of $\Omega$.
In many applications, e.g., optimal control and shape optimization,
the Dirichlet boundary data $g(\x)$ is rough, i.e., $g  \notin H^{{1}/{2}}(\Gamma)$.
This implies that the solution $u$ is not in $H^1(\Omega)$,
hence, it does not satisfy the standard variational formulation.
As a result, the conventional definition of the weak solution of \refe{pde}-\refe{pde-bc} must be modified.
The transposition method of Lions and Magenes \cite{Lions-Magenes} introduces
the very weak variational formulation: Seek $u \in L^2(\Omega)$, such that
\begin{align}
	\label{very-weak-solution}
	(u, \Delta v) = \langle g, \partial_\bn v \rangle - (f,v)\,, \quad  \forall v \in V \,,
\end{align}
where $V= H^2(\Omega) \cap H_0^1(\Omega)$.
It is easy to see that \refe{very-weak-solution} only requires the boundary data $g \in L^2(\Gamma)$,
as the test function $v$ is assumed to possess a higher regularity.
Due to its important applications,
numerical methods for solving the very weak solution to the elliptic/parabolic problems with rough boundary data
have been extensively studied, see \cite{Apel,Apel2-siam,Berggren2004,Bramble-King,French-King,French-King1993}. 
In particular,
the standard Lagrange finite element method(FEM) combined with $L^2(\Gamma)$-projection of the boundary data $g$
is widely used.
In the pioneering work \cite{Berggren2004}, Berggren rigorously analyzed this approach
and proved that the $u_h$ obtained by $L^2(\Gamma)$-projection converges to
the very weak solution $u$ for general Lipschitz polygonal/polyhedral domains.
However, the test space $V$ to define \refe{very-weak-solution} only applies to convex domains.
For a two dimensional nonconvex polygon, the important work of Apel, Nicaise and Pfefferer \cite{Apel} presents some remedies
by using an enlarged test space
$V = (H^2(\Omega) \cap H^1_0(\Omega)) \oplus \textrm{Span}\{\xi(r)r^{\lambda}\sin(\lambda \theta)  \}$,
where $(r, \theta)$ denotes the polar coordinate and $\lambda=\pi/\Theta$ with $\Theta$ being the re-entrant angle.
By using adaptive mesh strategy and the singular complement method,
they in \cite{Apel2-siam} improved the performance
of the Lagrange FEM for two dimensional nonconvex polygonal domains.
Moreover, Apel et al. in \cite{Apel} proposed a regularized approach
that introduces a sequence of regularized functions $\{g^h\} \in H^{{1}/{2}}(\Gamma)$
such that
$\lim_{h \rightarrow 0} \|g^h-g\|_{L^2(\Gamma)} = 0$.
Then, standard linear FEM can be applied with the boundary data $g^h$.
It should be noted that using Lagrange FEM to solve Poisson's equations with $L^2$ boundary
data needs to modify the original boundary data $g$, e.g., the $L^2(\Gamma)$-projection \cite{Berggren2004}.
Thus, an additional step is needed to preprocess the Dirichlet boundary data.
However, numerical evidences show that the $L^2(\Gamma)$-projection approach
may introduce certain artificial oscillations near the singular boundary points.
Similar approach can also be found in
numerical methods for elliptic problems with discontinuous Dirichlet boundary data,
see \cite{Cai-Yang2023}. For two dimensional elliptic equations with discontinuous
boundary data, Houston and Wihler in \cite{Houston-Wihler2012}
introduced a weak form in terms of weighted Sobolev spaces.
They proposed an interior penalty discontinuous Galerkin(DG) method,
where a posteriori error estimation is also derived.
It should be noted that the boundary data in the DG method is used implicitly
and regularization of $g$ is not needed.

In this paper, we use a Raviart--Thomas mixed FEM
to solve the problem \refe{pde}-\refe{pde-bc} with $L^2$ boundary data.
There have been extensive studies on Raviart--Thomas mixed FEMs,
see \cite{Boffi-Brezzi-Fortin,Gatica,Nedelec,Raviart_Thomas}.
Mixed FEMs have been widely used in boundary control
problems governed by elliptic PDEs, e.g., see \cite{ChenLiu2006,Garg-Porwal,Gong-Yan2011}.
For the model problem \refe{pde}-\refe{pde-bc},
the mixed method introduces an extra variable $\bs=\nabla u$.
Then, by using integration by parts,
there holds
\begin{align}
	(\bs, \bm{\chi})+( u , \D \, \bm{\chi}) = {\langle g, \bm{\chi}\cdot \mathbf{n} \rangle }
\end{align}
\color{black}
for any smooth functions $u$, $\bs$ and $\bm{\chi}$.
Therefore, the Dirichlet boundary data $g$ is used in an implicit way in the above weak formulation.
Motivated by this observation, we propose to use the lowest order mixed FEM
$\RT_{\! \! 0} \times \DG_{\! 0}$ to solve \refe{pde}-\refe{pde-bc}, see Subsection \ref{femmethod-sec2}.
However, all previous analyses of Raviart--Thomas mixed FEM
require that the boundary data $g \in H^{1/2}(\Gamma)$ at least.
In this work, we prove the numerical solution $u_h$ of  the lowest order mixed FEM
$\RT_{\! \! 0} \times \DG_{\! 0}$ converges to the very weak solution $u$
and establish an optimal error estimate.
The main difficulty in the analysis of the $\RT_{\! \! 0} \times \DG_{\! 0}$ mixed FEM
lies in the fact that the standard mixed variational form does not hold for $L^2(\Gamma)$ boundary data.
As the very weak solution $u \notin H^1(\Omega)$, one has $\bs \notin L^2(\Omega)$.
Consequently, the classical error estimate framework for mixed FEM cannot be applied.
The analysis of mixed FEM to \refe{pde}-\refe{pde-bc} is nonstandard.
In this work, we provide an optimal $L^2$-norm error estimate.
The analysis is based on a regularized approach.
We split the error into two parts: the regularization error $u-u^h$
and the approximation error $u^h - u_h$, where $u^h$ denotes the solution to the regularized problem.
Furthermore, the results can be extended to the case $g \in H^{s}(\Gamma)$ for $0 < s < 1/2$
with an improved convergence rate.

The rest of this paper is organized as follows.
In section \ref{femmethod}, we present a mixed FEM for solving the elliptic problem with rough
boundary data and main theoretical results.
In Section \ref{preparation}, we introduce some useful lemmas.
In Section \ref{proof-main-result}, we prove the optimal error estimate for the mixed FEM.
In Section \ref{more-smooth-data},  we extend the results to problems with more regular boundary data.
In Section \ref{numericalresults},
several numerical examples are provided to confirm our theoretical analysis and demonstrate
the effectiveness of the mixed FEM.
Some concluding remarks are given in section \ref{sec-conclusion}.

\section{A mixed FEM and main results}
\label{femmethod}
\setcounter{equation}{0}
We will introduce some standard notations and define
the very weak solution for general Lipschitz polygonal/polyhedral domains
in Subsection \ref{femmethod-sec1}.
Then, we present the mixed finite element method and main results on the convergence
in Subsection \ref{femmethod-sec2}.

\subsection{The very weak solution}
\label{femmethod-sec1}

We consider a bounded polygonal (for $d=2$) or polyhedral (for $d=3$) domain $\Omega \in \mathbb{R}^d,d=2,3$,
with a Lipschitz boundary $\Gamma$.
For any two functions $u$, $v \in {L}^{2}(\Omega)$,
we denote the ${L}^{2}(\Omega)$ inner product in domain $\Omega$ and the $L^2$-norm by
\begin{equation}
	(u,v)  = \,
	\int_{\Omega} u(\x) \, v(\x) \, {\mathrm{d}} \x ,
	\qquad {\left \| u \right \|_{L^2(\Omega)}} = (u,u)^{\frac{1}{2}} \,.
	\nn
\end{equation}
Similarly,  the inner product and norm on the boundary are defined by
\begin{equation}
	\langle g,\omega \rangle  = \,
	\int_{\Gamma} g(\x) \, \omega(\x) \, {\mathrm{d}} \x \,,
	\qquad {\left \| g \right \|_{{L}^2(\Gamma)}} = { \langle g,g \rangle }^{\frac{1}{2}} \,.
	\nn
\end{equation}
Let $W^{k,p}(\Omega)$ be the Sobolev space defined on $\Omega$,
and ${W}^{k,p}_0(\Omega)$ be the subspace of $W^{k,p}(\Omega)$
with zero trace. By conventional notations,
we define $H^{k}(\Omega) := W^{k,2}(\Omega)$
and ${H}^{k}_0(\Omega) := {W}^{k,2}_0(\Omega)$.
For a positive real number $s= k + w$, with $ w \in (0,1)$,
we define $H^s(\Omega)=(H^{k}(\Omega), H^{k+1}(\Omega))_{[w]}$
via the complex interpolation, see \textcolor{black}{\cite[Theorem 6.4.5]{Bergh}} and \cite{kellogg}.
To abbreviate notations, we use $\| \cdot \|_{L^2(\Omega)}$ and  $\| \cdot \|_{H^r(\Omega)}$
to denote the $L^2$ and $H^r$-norm of the inner product functions in the domain
$\Omega$, respectively.
Moreover, we define $\bH(\D,\Omega)$ by
\begin{align}
	& \bH(\D,\Omega):=
	\{ \bs| \bs \in [L^2(\Omega)]^d,~ \D \bs \in L^2(\Omega)  \} \,,
\end{align}
with norm $\| \bs\|_{\bH(\D)} = \| \bs\|_{\bH(\D,\Omega)}
:= (\|\bs\|^2_{L^2(\Omega)} + \|\D \bs\|^2_{L^2(\Omega)})^{\frac{1}{2}} $.
For simplicity, we omit $\Omega$ and define $\bH(\D) := \bH(\D,\Omega)$.

Now we introduce the definition of very weak solutions to the model problem
\refe{pde}-\refe{pde-bc} with $L^2$ boundary data $g$ for general polyhedral domains.
It should be noted that Apel, Nicaise and Pfefferer in \cite{Apel}
investigated the very weak solution for general two dimensional polygonal domains.
Based on the results in \cite{Apel}, we define the test space $V$ of very weak solutions by
\begin{equation}
  V = H_0^1(\Omega) \cap \{v\in L^2(\Omega) : \Delta v \in L^2(\Omega)\}
  \label{very-weak-test-space}.
\end{equation}
Since $\|v\|_{H^1(\Omega)} \le C\|\Delta v\|_{L^2(\Omega)}$ for $v\in V$,we can define the norm of $V$ as
\begin{equation}
  \|v\|_V =\|\Delta v\|_{L^2(\Omega)}.
\end{equation}

The following lemma addresses the very weak solution in general polygonal/polyhedral domains.

\begin{lemma}
	
\label{lemma-very-weak-solution}
Let $\Omega$ be a bounded Lipschitz polygonal or polyhedral domain.
Let $g \in L^2(\Gamma)$ and $f \in H^{-1}(\Omega)$,
then there exists a unique solution $u \in L^2(\Omega)$ satisfying
\begin{align}
  \label{very-weak-scheme}
  (u, \Delta v) = \langle g, \partial_\bn v \rangle - \textcolor{black}{(f,v)_{-1,1}}\,, \quad    \forall v \in V
  \,.
\end{align}
Moreover, there holds
\begin{equation}
	\|u\|_{L^2(\Omega)} \le C (\|g\|_{L^2(\Gamma)} +\|f\|_{H^{-1}(\Omega)}),
	\label{error-estimate-u}
\end{equation}
\textcolor{black}{
where $(\cdot,\cdot)_{-1,1} $ represents the duality pairing between $H^{-1}(\Omega)$ and $H_0^1(\Omega)$.}
\end{lemma}
{\noindent \bf Proof.}
Noting the regularity result in Lemma \ref{lemma-regularity} and
trace inequalities in Corollary \ref{trace-lemma}, we have
\begin{equation}
	\|v\|_{H^{1+s}(\Omega)} \le C\|\Delta v\|_{L^2(\Omega)} = C \| v\|_V .
\end{equation}
Hence, we obtain the embedding $V \hookrightarrow H^{s+1}(\Omega)\cap H_0^1(\Omega)$.
Since $s> 1/2$, combining Lemma \ref{trace-lemma} there holds
\begin{equation}
  \|\partial_n v\|_{L^2(\Gamma)} \leq C \|\partial_n v\|_{H^{s-\frac{1}{2}}(\Gamma)}
  \leq C \|v\|_{H^{1+s}(\Omega)} \leq C \|v\|_V,\qquad \forall v \in V.
\end{equation}
Then the right side of \ref{very-weak-scheme} defines a bounded linear functional on $V$.
To obtain the inf-sup condition, we can follow the same approach in \cite[Lemma 2.3]{Apel}, as
the proof remains valid in the three dimensional case.
\qed

\color{black}

It should be remarked that the above definition for the very weak solution is an extension of
the one defined in \cite{Apel} for two dimensional problems.
Due to the equivalence of $\|\Delta v\|_{L^2(\Omega)}$ and $\|v\|_{H^2(\Omega)}$ in convex domains, 
the test space $V$ in this case is the same to $H^2(\Omega) \cap H_0^1(\Omega)$.
\textcolor{black}{
In Berggren's work \cite{Berggren2004}, he also defined
the very weak solution for Poisson's equations on general Lipschitz polygonal or polyhedral domains.
If $g \in L^2(\Gamma)$ and $f \in L^2(\Omega)$,
Berggren's approach is equivalent to the above definition \refe{very-weak-scheme}.}
It should be noted that the source term $f \in H^{-1}(\Omega)$ may also introduce singularities.
However, the emphasis of this paper is on the error analysis of Poisson's problems with rough boundary data.
Thus, we shall assume $f \in L^2(\Omega)$.

\subsection{A mixed FEM and main results on error estimate}
\label{femmethod-sec2}

Let $\mathcal{T}_{h}$ be a quasi-uniform tetrahedral partition(triangular partition in 2D) of $\Omega$
with $\Omega = \cup_{K  \in \mathcal{T}_{h} } \Omega_{K}$ and denote by
$h=\max_{\Omega_{K} \in \mathcal{T}_{h}}\{ \mathrm{diam} \, \Omega_{K}\}$ the
mesh size.
By $\mathcal{F}_h$ we denote all the $(d-1)$-dimensional faces of the mesh
partition $\mathcal{T}_h$. Let $\mathcal{F}_h^{\partial}=\mathcal{F}_h \cap \Gamma$.
For $r \ge 0$, we define the Raviart--Thomas mixed
finite element spaces by
\begin{empheq}[left=\empheqlbrace]{align}
	& \RT_{\! \! r} := \{ \bm{\chi}_h \in \bH(\D): \bm{\chi}_h |_K \in [P_r(K)]^d + \x  P_r(K), \forall K
	\in \mathcal{T}_{h} \}\,,
	\nn\\
	& {\DG}_r := \{ p_h \in  L^2(\Omega): p_h |_K \in P_r(K), \forall K \in \mathcal{T}_{h} \}\,,
	\nn
\end{empheq}
where $P_r(K)$ is the space of polynomials of degree $r$ or less defined on $K$.
It is well-known that $\RT_{\! \! r} \times {\DG}_r $ is a stable finite
element pair
for second order elliptic problems, see \cite{Boffi-Brezzi-Fortin,Gatica,Nedelec,Raviart_Thomas}.

With the above notations, a mixed FEM for \refe{pde}-\refe{pde-bc}
is to seek $(\bs_h,u_h) \in \RT_{\! \! 0} \times \DG_{\! 0} $, such that
\begin{empheq}[left=\empheqlbrace]{align}
	&(\bs_h ,\, \bm{\chi}_h) + (u_h  , \, \D \bm{\chi}_h )= \langle g ,\, \bm{\chi}_h
	\cdot \n \rangle \,,
	&& \forall \bm{\chi}_h \in \RT_{\! \! 0}\,, 
	\label{mixed-fem1}
	\\
	& -(\D \bs_h  ,\, v_h) = (f,v_h)\,,
	&&
	\forall \, v_h \in \DG_{\! 0}\,.
	\label{mixed-fem2}
\end{empheq}
By noting the fact that $\bm{\chi}_h \cdot \n$ is piecewise constant on $\mathcal{F}_h^{\partial}$,
the inner product $\langle g, \bm{\chi}_h \cdot \n \rangle$ is
well-defined for any $g \in L^{2}(\Gamma)$.
The existence and uniqueness of the numerical solution $(\bs_h,u_h)$ to
\refe{mixed-fem1}-\refe{mixed-fem2} have been well studied, see
\cite{Boffi-Brezzi-Fortin}.
In addition, it should be pointed out that higher order elements are not
useful as the exact solution $u \notin H^1(\Omega)$.

We present our main results for the mixed FEM \refe{mixed-fem1}-\refe{mixed-fem2} in the following
theorem. The proof will be given in Section \ref{proof-main-result}.

\begin{theorem}
	\label{main-theorem}
	
	Let $f \in L^2(\Omega)$, $g \in L^2(\Gamma)$,
	the mixed FEM \refe{mixed-fem1}-\refe{mixed-fem2} admits a unique solution $u_h$ which
	converges to the very weak solution $u$ defined in \refe{very-weak-scheme},and there holds
	\begin{align}
		\| u_{h} - u \|_{L^2(\Omega)} \le C h^{s-\frac{1}{2}}\|g\|_{L^2(\Gamma)} + C h^s \|f\|_{L^2(\Omega)}\,,
		\label{main-error-estimate-u}
	\end{align}
	where $C$ is a positive constant independent of $h$,
	and the index $s>1/2$ is defined in \refe{index-unified}.
\end{theorem}

In the rest of this paper, we denote by $C$ a generic positive constant
and by $\ee$ a generic small positive constant, which are independent of $h$.

\vspace{0.4cm}

%%%%%%%%%%%% 3.1
\section{Preliminaries}
\label{preparation}
\setcounter{equation}{0}

In this section, we present several useful lemmas, which will be frequently used in our proof.
Let $\P_h : L^2(\Omega) \rightarrow {\DG}_0$ be the $L^2$ projector: For $u\in L^2(\Omega)$, seek $\P_h u \in {\DG}_0$, such that
\begin{align}
	(\mathcal{P}_h u - u , v_h) =0,  \quad \forall ~ v_h \in {\DG}_0 \,.\label{equ-dg0-projector}
\end{align}
Let ${\Pi}_h: \bH(\D) \rightarrow {\RT_{\! \! 0} }$ be the quasi Raviart--Thomas projector
developed by Ern et al. in \cite{Ern-Gudi}.
Then, the following diagram commutes \cite{Ern-Gudi}
\begin{align}
	\begin{array}{crlrlrl}%ccccccc
		\bH(\D)
		& \xrightarrow{\quad \D \quad}
		& \quad L^2(\Omega)              \\
		\Pi_h~ \!\!\Bigg\downarrow ~~~~ &                              &
		~ \P_h~ \!\! \Bigg\downarrow                                     \\
		\RT_{\! \! 0}
		& \xrightarrow{\quad \D \quad}
		& \quad~ {\DG}_0
	\end{array}
	\label{diagram}
\end{align}
Moreover, the following lemma holds under the minimal necessary Sobolev regularity \cite{Ern-Gudi}.
\begin{lemma}
	The quasi-projector $\Pi_h$ maps $\bH(\D)$ to $\RT_{\! \! 0}$ and there holds
	\begin{align}
		\label{RTNp-error}
		\| \bv - \Pi_h\bv\|^2_{L^2(\Omega)}
		& +  h^2 \| \D \, (\bv - \Pi_h\bv )\|_{L^2(\Omega)}^2
		\nn                                                       \\
		& \le C \left( h^{\min{(s,1)}} \|\bv\|_{\bH^s(\Omega)} +
		\delta_{s<1} \, h \, \|\nabla\cdot\bv\|_{L^2(\Omega)} \right)^2\,,
	\end{align}
	where $\delta_{s<1}:=1$ if $s<1$ and $\delta_{s<1}:=0$ if $s\ge1$.
	In addition, the projector $\Pi_h$ is globally $L^2$-stable up to $hp$
	data oscillation of the divergence and $\bH(\D)$-stable
	\begin{empheq}[left=\empheqlbrace]{align}
		& \|\Pi_h\bv\|^2_{L^2(\Omega)} \le C \left( \|\bv\|^2_{L^2(\Omega)}
		+ h^2 \| \D \, \bv - \P_h (\D \, \bv) \|_{L^2(\Omega)}^2 \right) \, ,
		\label{RTNp-stable}
		\\
		&
		\|\Pi_h \bv\|^2_{L^2(\Omega)} +  \|\D \, \Pi_h \bv\|^2_{L^2(\Omega)}
		\le C \left( \|\bv\|^2_{L^2(\Omega)} + \|\D \, \bv\|^2_{L^2(\Omega)} \right)\,.
		\label{RTNp-stable-Hdiv}
	\end{empheq}
	
\end{lemma}

Moreover, the following error estimates hold for $\Pi_h$ and $\P_h$, see
\cite{Ciarlet} and \cite[Section 3]{ACWY}
\begin{empheq}[left=\empheqlbrace]{align}
	& \| v- \P_h v\|_{H^{-t}(\Omega)} \le C h^{t+s} \|v\|_{H^s(\Omega)} \,,
	&& \textrm{for}~ 0 \le t,s \le 1\,.
	\label{L2-projection-error-u}
	\\
	& \|\D(\bm{v} -\Pi_h \bm{v}) \|_{L^2(\Omega)} \le  C h^{s} \|\D \, \bm{v}\|_{H^s(\Omega)}\,,
	&& \textrm{for}~ 0 \le s \le 1  \,,
	\\
	& \| (\bm{v} - \Pi_h \bm{v}) \cdot  \bn\|_{H^{-t}(\Gamma)} \le Ch^{t+s}  \|\bm{v}\|_{H^s(\Gamma)}\,,
	&& \textrm{for}~ 0 \le t,s \le 1 \, ,
	\label{RT-projection-error-boundary}
\end{empheq}
The below inverse estimate for the normal trace holds \cite[Lemma 4.1]{ACWY}
\begin{align}
	\| \bm{v}_h \cdot \bn \|_{L^2(\Gamma)} \le C h^{-\frac{1}{2}}
	\|\bm{v}_h\|_{L^2(\Omega)}\,,
	\qquad
	\forall \bm{v}_h \in \RT_{\! \! 0}  \,.
	\label{RT-element-boundness-boundary}
\end{align}

The following results on traces are needed in our analysis,
see \cite[Theorem 1.5.1.2,Theorem 1.5.1.3]{Grisvard1985EllipticPI}.

\begin{lemma}
	\label{trace-lemma-origin}
	Let $\Omega$ be a bounded Lipschitz domain.
	Assume $\mathrm{tr}:H^s(\Omega) \to L^2_{loc}(\Gamma)$ is the trace operator on $\Gamma$,
	then for $u\in H^s(\Omega)$ with ${1}/{2} < s\le 1$, there holds
	\begin{equation}
		\|\mathrm{tr}(u)\|_{H^{s-\frac{1}{2}}(\Gamma)} \le C \|u\|_{H^{s}(\Omega)}\,,
	\end{equation}
	where C depends on the domain $\Omega$ only.
\end{lemma}

\begin{corollary}
	\label{trace-lemma}
	Let $\Omega$ be a bounded Lipschitz domain.  
	For $u\in H^{1+s}(\Omega)$ with ${1}/{2} < s \le 1$, there holds
	\begin{equation}
		\|\mathrm{tr}(\frac{\partial u}{\partial \bn})\|_{H^{s-\frac{1}{2}}(\Gamma)} \le C \|u\|_{H^{1+s}(\Omega)}\,,
	\end{equation}
	where $C$ depends on $\Omega$ only, $\bn$ is the outer normal on the boundary $\Gamma$.
\end{corollary}

\begin{lemma}
  \label{lemma-regularity}
  The solution $u$ to the Poisson equation with a homogeneous Dirichlet boundary condition
  \begin{equation}
    \left\{ {
      \begin{array}{l}
	{-\Delta u =  f \,, \quad \quad \mathrm{in} ~ \Omega\,, } \\
	{u = 0\,, \qquad\quad  ~\mathrm{on}~ \Gamma\,,}
    \end{array}}
    \right.
    \label{poisson-standard}
  \end{equation}
  satisfies
  \begin{align}
    \|u\|_{H^{1+\lambda}(\Omega)} \le C \|f\|_{L^2(\Omega)},
  \end{align}
  where $f\in L^2(\Omega)$, $C$ is a positive constant independent of $u$ and $f$.
  The index $\lambda$ only depends on the domain $\Omega$.
  For two dimensional polygons, we have
  \begin{equation}
    \lambda_{\textrm{2D}} \in \big( \frac{1}{2}, \, \frac{\pi}{\max_j \Theta_j} \big)  \, ,
    \label{s-index-2d}
  \end{equation}
  where $\{\Theta_j\}$ denotes the re-entrant interior angles of $\Omega$.
  For three dimensional nonconvex polyhedral domain, the regularity of $u$ depends on both
  the edge opening angle at edges $\{ e \}$ and the shape of the domain near corners $\{v\}$,
  i.e., edges and corners may introduce certain singularities.
  \textcolor{black}
  {Assume $\mathcal{E}$ and $\mathcal{V}$ represent all the edges and vertices respectively,
  the solution $u$ satisfies the splitting
  \begin{equation}
  	u = u_r + \sum_{e\in \mathcal{E}} \alpha_e \psi_e u_e +  \sum_{v\in \mathcal{V}} \alpha_v \psi_v u_v \,,
  \end{equation}
  where $u_r \in H^2(\Omega)$ denotes the regular part,
  $\psi_e$ and $\psi_v$ are cutoff functions
  that equal 1 in neighborhoods of $e$ and $v$, respectively. 
  Here, $u_e$ and $u_v$ denote the singular functions associated with the edge and vertex.
  $\alpha_e$ and $\alpha_v$ represent their corresponding singularity coefficients.}
  \color{black}
  Moreover, there holds
  \begin{equation}
    \label{s-index-3d-part1}
    u_e \in H^{1+\lambda_{\textrm{3D}}}(\Omega),
    \quad \textrm{with} \quad
    \frac{1}{2} < \lambda_{\textrm{3D}} < \frac{\pi}{\max_{e \in \mathcal{E}} \Theta_e} \,,
  \end{equation}
  and 
  \begin{equation}
    \label{s-index-3d-part2}
    \begin{aligned}
    u_v \in H^{1+\lambda_{\textrm{3D}}}(\Omega),
    \quad \textrm{with} \quad      
      & \frac{1}{2} < \lambda_{\textrm{3D}} < \frac{1}{2} + \min_{v \in \mathcal{V}}\{ \lambda_{v,D} , 2\},
    \end{aligned}
  \end{equation}
  \color{black}
  where $\Theta_e$ denotes the edge opening angle at the edge $e$,
  and $\lambda_{v,D}>0$ depends on the smallest Dirichlet eigenvalue of the surface Laplacian
  around corner $v$, see \cite{Chatzipantelidis-BIT2006,Dauge,Grisvard1992Singular} for details.
  
\end{lemma}

If the domain $\Omega$ is convex, there holds
\begin{align}
	\|u\|_{H^{2}(\Omega)} \le C \|f\|_{L^2(\Omega)}.
\end{align}

In the rest of this paper, we use a unified index $s$ to describe the regularity of \refe{poisson-standard}
\begin{align}
  \|u\|_{H^{1+s}(\Omega)} \le C \|f\|_{L^2(\Omega)}, \quad
  \textrm{with} \quad
  \left\{
  \begin{array}{ll}
  	s = 1
  	& \textrm{if $\Omega$ is convex,}
  	\\[2pt]
  	s = \sup \{ \lambda_{\textrm{2D}} \} - \epsilon
  	& \textrm{in 2D},
  	\\[2pt]
  	s = \sup \{ \lambda_{\textrm{3D}} \} - \epsilon
  	& \textrm{in 3D},
  \end{array}
  \right.
  \label{index-unified}
\end{align}
where
\textcolor{black}{
  $\lambda_{\textrm{2D}}$ is defined in \refe{s-index-2d}
  and
  $\lambda_{\textrm{3D}}$ is defined in \refe{s-index-3d-part1}-\refe{s-index-3d-part2},
  respectively,}
and  $\epsilon>0$ is any arbitrarily small number.
We can see that for Lipschitz polygonal/polyhedral domains,
the solution $u$ to \refe{poisson-standard} always belongs to $H^{3/2+\epsilon}(\Omega)$,
\textcolor{black}{see \cite[Theorem 3.1]{Berggren2004}}.
Recall the definition of $V$ in \refe{very-weak-test-space},
the $H^{3/2+\epsilon}$ regularity 
ensures that $\partial_{\mathbf{n}} v \in L^{2}(\Gamma)$ for any $v \in V$.

The standard Raviart--Thomas mixed FEM for \refe{poisson-standard} is to
seek $(\bs_h,u_h) \in \RT_{\! \! 0} \times \DG_{\! 0} $, such that
\begin{empheq}[left=\empheqlbrace]{align}
	&(\bs_h \,,\, \bm{\chi}_h) + (u_h \,, \D \bm{\chi}_h )= 0 \,,
	&& \forall \bm{\chi}_h \in \RT_{\! \! 0} \,,
	\label{mixed-fem1-standard}
	\\
	& -(\D \bs_h, v_h) = (f,v_h),
	&&
	\forall \,v_h \in \DG_{\! 0} \, .
	\label{mixed-fem2-standard}
\end{empheq}
The following estimates hold \cite[Lemma 6.1]{Ern-Gudi}
\begin{align}
	\|\bs_h - \bs\|_{L^2(\Omega)} \le C h^s  \|f\|_{L^2(\Omega)}  \quad
	\textrm{and} \quad
	\|u_h - u \|_{L^2(\Omega)} \le C h \|f\|_{L^2(\Omega)}\,,
\end{align}
where the index $s$ is defined in \refe{index-unified}.

\vspace{0.4cm}

%==============================
% The proof of error result
%==============================
\section{The proof of the main result}
\label{proof-main-result}
\setcounter{equation}{0}

Following the idea of Apel, Nicaise and Pfefferer in \cite{Apel}, we introduce a regularized
elliptic problem for the original Poisson equation \refe{pde}-\refe{pde-bc}.
We shall introduce the linear Lagrange element space $P_1$ on the mesh $\mathcal{T}_{h}$.
Moreover, let $P_1^\partial:=P_1|_{\Gamma}$.
Let $g^h \in H^{1/2}(\Gamma)$ denote a sequence of functions such that
\begin{align}
	\lim_{h \to 0} \|g - g^h\|_{L^2(\Gamma)} = 0.
\end{align}
The construction of $g^h$ can be done by the $L^2$-projection of $g$ onto $P_1^\partial$,
which is analyzed in Berggren's pioneering work \cite{Berggren2004}.
Alternatively, one can use the Carstensen interpolant \cite{Carstensen1999}.
\color{black}
For the sequence $\{g^h\}$, there holds \cite{Berggren2004},\cite[Lemma 2.14]{Apel}\color{black}
\begin{align}
	\label{gh-error}
	\| g^h - g \|_{H^{-s}(\Gamma)} \le Ch^{t+s}\|g\|_{H^t(\Gamma)}, \quad
	\forall \, t, ~s \in [0,1],
\end{align}
provided $g \in H^t(\Gamma)$.

Now, for given $g^h \in H^{{1}/{2}}(\Gamma)$, $f \in L^2(\Omega)$,
we introduce a mixed weak form: Seek $\bs^h \in \bH(\D), u^h \in L^2(\Omega)$, such that
\begin{empheq}[left=\empheqlbrace]{align}
	&( \bs^h ,\,  \bm{\chi})
	+(u^h,\, \D  \bm{\chi})= \langle g^h, \, \bm{\chi} \cdot \n \rangle \, ,
	&& \forall ~ \bm{\chi} \in \bH(\D)\,,
	\label{auxiliary-variational-problem-1}
	\\
	&(-\D \bs^h,\, v) = (f, \, v) \,,
	&& \forall ~ v \in L^2(\Omega)\,.
	\label{auxiliary-variational-problem-2}
\end{empheq}
The existence of the solution $(\bs^h,u^h)$ to the above mixed variational form
is obvious \cite{Boffi-Brezzi-Fortin,Gatica}, which satisfies the following standard estimate
\begin{align}
	\|\bs^h\|_{\bH(\D)} + \|u^h\|_{L^2(\Omega)} \le C\left( \|f\|_{L^2(\Omega)}+ \|g^h\|_{H^{\frac{1}{2}}(\Gamma)}\right).
	\label{mixedFEM-standard-estimate}
\end{align}
The following lemma shows the convergence of solutions $\{u^h\}$
of the auxiliary problem to the very weak solution $u$.	
\color{black}
	This result follows directly from \cite[Corollary 3.3]{Apel}, since the unique weak solution
	to the primal formulation is also the unique weak solution to its corresponding mixed formulation.
	While the original corollary is stated for polygonal domains,
	the proof is still valid for Lipschitz polyhedral domain
	as Lemma \ref{trace-lemma} provides the required trace regularity.
\color{black}
\begin{lemma}
	\label{error-auxiliary-variational-problem}
	For given $g^h \in \{ g^h \}$, the variational form
	\refe{auxiliary-variational-problem-1}-\refe{auxiliary-variational-problem-2}
	admits a unique solution $u^h \in L^2(\Omega)$.
	The limit $u := \lim_{h \rightarrow 0} u^h$  exists and is the very weak solution.
	Moreover, there holds,
	\begin{align}
		\label{error-auxiliary-uh}
		\| u^h - u \|_{L^2(\Omega)}  \le C h^{s-\frac{1}{2}}  \|g\|_{L^2(\Gamma)} \,,
	\end{align}
	where the index $s$ is defined in \refe{index-unified}.  
\end{lemma}

Based on the regularization error \refe{error-auxiliary-uh}
in Lemma \ref{error-auxiliary-variational-problem},
we prove the error estimate \refe{main-error-estimate-u}
in the main Theorem \ref{main-theorem}.
The proof consists of two steps. In the first step,
we derive an estimate for $\bs^h-\bs_h$.
In the second part, we present an error estimate for $u^h-u_h$.

By subtracting the regularized equations
\refe{auxiliary-variational-problem-1}-\refe{auxiliary-variational-problem-2}
from mixed FEM \refe{mixed-fem1}-\refe{mixed-fem2}, we deduce the error equations
\begin{empheq}[left=\empheqlbrace]{align}
	\label{error-pre-1}
	&(\bs^h - \bs_h , \bm{\chi}_h) + (u^h-u_h , \D \bm{\chi}_h)
	= \langle g^h-g\,, \bm{\chi}_h \cdot \bn \rangle ,
	&& \forall \bm{\chi}_h \in \RT_{\! \! 0} \,,
	\\
	\label{error-pre-2}
	& (\D (\bs^h-\bs_h) , v_h) = 0 ,
	&& \forall v_h \in \DG_0 \,.
\end{empheq}
Here, we shall introduce the $L^2$ projector of $u^h$,
\begin{align}
	(\P_h u^h - u^h , v_h ) = 0,   \qquad \forall  v_h \in \DG_{\! 0} \,.
	\nn
\end{align}
Then, the projection error satisfies
\begin{eqnarray}
	\|\P_h u^h - u^h\|_{L^2(\Omega)} & \le & C h \|u^h\|_{H^1(\Omega)}
	\nn \\
	& \le & C h (\|f\|_{L^2(\Omega)} + \|g^h\|_{H^{\frac{1}{2}}(\Gamma)})
	\nn \\
	& \le & C h \|f\|_{L^2(\Omega)} + C h^{\frac{1}{2}}\|g\|_{L^{2}(\Gamma)}.
	\label{error-uh-mid}
\end{eqnarray}
where we have used an inverse inequality and the fact that $g^h \in P_1^{\partial}$.
As $\D \RT_{\! \! 0} \subset \DG_{\! 0} $, the error equations
\refe{error-pre-1}-\refe{error-pre-2} can be rewritten as
\begin{empheq}[left=\empheqlbrace]{align}
	& (\bs^h - \bs_h , \bm{\chi}_h)  +(\P_h u^h -u_h , \D \bm{\chi}_h)
	= \langle g^h-g\,, \bm{\chi}_h  \cdot \bn \rangle ,  && \forall  \bm{\chi}_h \in \RT_{\! \! 0} \,,
	\label{error1}
	\\
	& (\D (\bs^h-\bs_h) , v_h) = 0 ,
	&&   \forall v_h \in \DG_0 \,.
	\label{error2}
\end{empheq}

An estimate for $\bs_h$ is given in the lemma below.
\begin{lemma}
	\label{lemma-sigma-eror}
	If $f \in L^2(\Omega)$, $g \in
	L^2(\Gamma)$ and $g^h \in H^{{1}/{2}}(\Gamma)$, we have
	\begin{align}
		\label{error-sigma}
		\| \bs^h - \bs_h\|_{L^2(\Omega)} \le C( h^{-\frac{1}{2}} \|g\|_{L^2(\Gamma)} + \|f\|_{L^2(\Omega)})\,,
	\end{align}
	where the constant $C$ is independent of $h$.
\end{lemma}
\
{\noindent \bf Proof.}
From the commuting diagram \refe{diagram}, the equation \refe{error2} can be rewritten as
\begin{eqnarray}
	(\D \, (\bs^h - \bs_h), v_h)
	& = &
	(\D \, (\bs^h -  \Pi_h \bs^h), v_h )
	+(\D \, (\Pi_h\bs^h -  \bs_h), v_h)
	\nn \\
	& = & (\D (\Pi_h\bs^h -  \bs_h), v_h)
	=0,
	\quad \forall v_h \in \DG_{\! 0}\,,
	\nn
\end{eqnarray}
which, by taking $v_h = \D ( \Pi_h\bs^h -\bs_h ) $ implies the fact
\begin{align}
	\label{div-sigma-free}
	\D (\Pi_h\bs^h - \bs_h)=0\,.
\end{align}
Next, by taking $\bm{\chi}_h = \Pi_h \bs^h - \bs_h$ into \refe{error1}, we obtain
\begin{align}
	\|\Pi_h\bs^h - \bs_h\|^2_{L^2(\Omega)}
	& =(\Pi_h\bs^h - \bs^h,  \Pi_h\bs^h - \bs_h)
	+ \langle g^h - g, (\Pi_h\bs^h - \bs_h) \cdot \bn \rangle
	\nn \\
	& \le \|\Pi_h\bs^h - \bs^h\|_{L^2(\Omega)}
	\|\Pi_h\bs^h - \bs_h\|_{L^2(\Omega)}
	\nn \\
	&  \quad
	+ \| g^h - g\|_{L^2(\Gamma)} \|(\Pi_h\bs^h-\bs_h) \cdot
	\bn\|_{L^2(\Gamma)}
	\nn \\
	\textrm{(by \refe{RTNp-error})} \quad
	& \le C (\|\bs^h\|_{L^2(\Omega)} + h \|\D \bs^h\|_{L^2(\Omega)} )
	\|\Pi_h\bs^h -  \bs_h\|_{L^2(\Omega)}
	\nn \\
	\textrm{(by inverse inquality)} \quad
	& \quad + Ch^{ -\frac{1}{2} } \| g^h - g\|_{L^2(\Gamma)}
	\|(\Pi_h\bs^h-\bs_h) \cdot
	\bn\|_{H^{-\frac{1}{2}}(\Gamma)}
	\nn \\
	& \le C \|\bs^h\|_{\bH(\D)}
	\|\Pi_h\bs^h - \bs_h\|_{L^2(\Omega)}
	+ Ch^{-\frac{1}{2}} \|g\|_{L^2(\Gamma)} \| \Pi_h\bs^h - \bs_h\|_{\bH(\D)}
	\nn \\
	\textrm{(by \refe{mixedFEM-standard-estimate})} \quad
	& \le C( h^{-\frac{1}{2}}  \|g\|_{L^2(\Gamma)} + \|f\|_{L^2(\Omega)} +\|g^h\|_{H^{\frac{1}{2}}(\Gamma)})
	\| \Pi_h\bs^h - \bs_h\|_{L^2(\Omega)}\,,
	\nn
\end{align}
which, by using an inverse inequality for $g^h$, leads to the result below
\begin{align}
	\|\Pi_h\bs^h - \bs_h\|_{L^2(\Omega)}  \le C( h^{-\frac{1}{2}}  \|g\|_{L^2(\Gamma)} + \|f\|_{L^2(\Omega)}).
	\label{final-sigma-error}
\end{align}
From \refe{RTNp-error}, we see that the quasi-projection estimate holds
\begin{eqnarray}
	\|\Pi_h\bs^h - \bs^h\|_{L^2(\Omega)} & \le & C\|\bs^h\|_{L^2(\Omega)}+Ch\|\D \, \bs^h\|_{L^2(\Omega)}
	\nn
	\\
	& \le & C (\|g^h\|_{H^{\frac{1}{2}}(\Gamma)}+\|f\|_{L^2(\Omega)})+Ch\|f\|_{L^2(\Omega)}
	\nn
	\\
	& \le & C h^{-\frac{1}{2}}\|g\|_{L^2(\Gamma)} + C\|f\|_{L^2(\Omega)}.
	\label{quasi-projection-sigma-error}
\end{eqnarray}

By combining estimate \refe{final-sigma-error}
and projection estimate \refe{quasi-projection-sigma-error},
the error estimate \refe{error-sigma} is proved.
\qed

Next, we turn to prove the main error estimate \refe{main-error-estimate-u} for $u_h$.

{\noindent \bf Proof.}
Notice that
\begin{align}
	\|u-u_h\|_{L^2(\Omega)} \le \|u-u^h\|_{L^2(\Omega)}  + \|u^h - \P_h u^h\|_{L^2(\Omega)}
	+ \|\P_h u^h - u_h\|_{L^2(\Omega)}\,.
\end{align}
Clearly, we only need to estimate the last term.
To this end, we introduce a dual Poisson's equation with a homogeneous Dirichlet boundary condition,
\begin{empheq}[left=\empheqlbrace]{align}
	& - \Delta z = \P_h u^h - u_h \,, && \textrm{in} ~ \Omega\,,
	\\ \nn
	& z =0\,,  && \textrm{on}~  \partial \Omega \,.
\end{empheq}
By introducing  $\bm{\omega}= \nabla z$ as an extra variable,
the mixed FEM for the above equation is to seek $(\bm{\omega}_h, z_h) \in \RT_{\! \! 0} \times \DG_{\! 0}$
such that
\begin{empheq}[left=\empheqlbrace]{align}
	& (\bm{\omega}_h, \bm{\nu}_h ) + (z_h, \D \bm{\nu}_h ) = 0 \,,
	&& \forall \bm{\nu}_h \in \RT_{\! \! 0}\,,
	\\
	& -(\D \bm{\omega}_h,  v_h ) = (\P_h u^h - u_h, v_h)\,,
	&& \forall v_h \in \DG_{\! 0} \,.
	\label{error-standard-proof}
\end{empheq}
By the standard error estimate of mixed FEMs for elliptic equation \cite{Ern-Gudi}, there holds
\begin{empheq}[left=\empheqlbrace]{align}
	\label{standard-mixed-fems-error}
	& \|\nabla z -\bm{\omega}_h\|_{L^2(\Omega)} \le
	C \left( h^s \|\bm{\omega}\|_{\bH^s(\Omega)} + h \|\D \bm{\omega}\|_{L^2(\Omega)} \right) \le C h^s \|\P_h u - u_h\|_{L^2(\Omega)} \,,
	\nn \\
	& \| \Pi_h \nabla z - \bm{\omega}_h \|_{L^2(\Omega)} \le C h^s \, \| \P_h u^h - u_h \|_{L^2(\Omega)} \,,
	\\
	& \| z-z_h \|_{L^2(\Omega)} \le C h\| \P_h u - u_h\|_{L^2(\Omega)}\,.
	\nn
\end{empheq}
By taking $v_h=\P_h u^h - u_h$ into \refe{error-standard-proof}, we can see that
the $L^2$-norm of $\P_h u^h - u_h$ satisfies
\begin{align}
	\label{error-u}
	\|\P_h u^h - u_h\|_{L^2(\Omega)}^2
	& = -(\D \bm{\omega}_h ,\P_h u^h - u_h)
	\nn                                                                                       \\
	\textrm{(by \refe{error1})} \quad
	& =  \langle g-g^h, \bm{\omega}_{h} \cdot \bn \rangle  + (\bs^h - \bs_h , \bm{\omega}_h)
	\nn                                                                                       \\
	& = \mathcal{J}_1 + \mathcal{J}_2\,.
\end{align}
By using the trace inequality in Corollary \ref{trace-lemma},
the projector error for the normal trace \refe{RT-projection-error-boundary}
and the inverse inequality \refe{RT-element-boundness-boundary},
the term $\mathcal{J}_1$ can be bounded by
\begin{align}
	\label{error-J1}
	\left| \mathcal{J}_1 \right|
	& = \left| \langle g-g^h, (\bm{\omega}_{h}- \Pi_h \nabla z) \cdot \bn \rangle
	+ \langle g-g^h, (\Pi_h \nabla z - \nabla z) \cdot \bn \rangle
	+ \langle g-g^h, \nabla z \cdot \bn  \rangle  \right|
	\nn
	\\[3.5pt]
	& \le \|g-g^h\|_{L^2(\Gamma)} \|(\bm{\omega}_{h}- \Pi_h \nabla z) \cdot \bn \|_{L^2(\Gamma)}
	+ \|g-g^h\|_{L^2(\Gamma)} \|( \Pi_h \nabla z - \nabla z )
	\cdot \bn \|_{L^2(\Gamma)}
	\nn
	\\[3.5pt]
	& \quad + \|g-g^h\|_{H^{-(s-\frac{1}{2})}(\Gamma)} \|\nabla z \cdot \bn \|_{H^{(s-\frac{1}{2})}(\Gamma)}
	\nn \\
	& \le C h^{-\frac{1}{2}}\|(\bm{\omega}_{h}- \Pi_h \nabla z) \|_{L^2(\Omega)} \|g\|_{L^2(\Gamma)}
	+ C h^{s-\frac{1}{2}} \|\nabla z\|_{H^{s-\frac{1}{2}}(\Gamma)} \|g\|_{L^2(\Gamma)}
	\nn \\
	&  \quad
	+ C h^{s-\frac{1}{2}} \|g\|_{L^2(\Gamma)} \| z\|_{H^{s+1}(\Omega)}
	\nn \\
	& \le   C h^{s-\frac{1}{2}} \|\P_h u^h  -u_h\|_{L^2(\Omega)}  \|g\|_{L^2(\Gamma)}
	+ C h^{s-\frac{1}{2}} \|\nabla z\|_{H^{s}(\Omega)} \|g\|_{L^2(\Gamma)}
	\nn \\
	&  \quad
	+ C h^{s-\frac{1}{2}} \|g\|_{L^2(\Gamma)} \| z\|_{H^{s+1}(\Omega)}
	\nn \\
	& \le C h^{s-\frac{1}{2}} \|g\|_{L^2(\Gamma)}
	\|\P_h u^h -u_h\|_{L^2(\Omega)} \, .
\end{align}
By the standard error estimate \refe{standard-mixed-fems-error}, the term $\mathcal{J}_2$
can be bounded by
\begin{align}
	\label{error-J2}
	\left| \mathcal{J}_2 \right|
	& \le |(\bs^h  -\bs_h , \bm{\omega}_h - \nabla z)| + |(\bs^h-\bs_h , \nabla z  )|
	\nn \\
	\textrm{(by \textrm{\refe{error2}})} \quad
	& = |(\bs^h  -\bs_h , \bm{\omega}_h - \nabla z) | + |(\D(\bs^h - \bs_h) , z - \P_h z)|
	\nn \\[3.5pt]
	& \le \|\bs^h-\bs_h\|_{L^2(\Omega)} \|\bm{\omega}_h - \nabla z\|_{L^2(\Omega)}
	+ \|\D(\bs^h - \bs_h)\|_{L^2(\Omega)} \| z - \P_h z\|_{L^2(\Omega)}
	\nn \\[3.5pt]
	\textrm{(by \refe{error-sigma} and \refe{standard-mixed-fems-error})} \quad
	& \le  C h^s  \left(h^{-\frac{1}{2}}\|g\|_{L^2(\Gamma)} + C\|f\|_{L^2(\Omega)} \right)
	\|\P_h u^h - u_h\|_{L^2(\Omega)}
	\nn \\
	& \quad
	+ Ch \|f\|_{L^2(\Omega)} \|z\|_{H^{1}(\Omega) }
	\nn \\
	& \le C\left(h^s\|f\|_{L^2(\Omega)} + h^{s-\frac{1}{2}}\|g\|_{L^2(\Gamma)}   \right)
	\|\P_h u^h - u_h\|_{L^2(\Omega)} \,.
\end{align}
Taking estimates \refe{error-J1} and \refe{error-J2} into \refe{error-u},
the desired estimate follows
\begin{align}
	\label{error-uh-pre}
	\| \P_h u^h - u_h \|_{L^2(\Omega)}
	\le C h^{s- \frac{1}{2}} \|g\|_{L^2(\Gamma)} + C h^s \|f\|_{L^2(\Omega)} \,.
\end{align}

Finally, combining estimates in \refe{error-estimate-u},
\refe{error-auxiliary-uh}, \refe{error-uh-mid} and
\refe{error-uh-pre}, the main results in Theorem \ref{main-theorem} is proved.
\qquad \endproof

\vspace{0.4cm}

%==============================
% Section on better g
%==============================
% \section{Estimates for problems with more regular boundary data}
% \label{more-smooth-data}
% \setcounter{equation}{0}

% Classical mixed FEM theory covers the case $g \in H^{t}(\Gamma)$ with $t \ge 1/2$.
% In this section, we investigate the model problem \refe{pde}-\refe{pde-bc}
% with boundary data $g \in H^{t}(\Gamma)$ for $0 < t < 1/2$.
% An improved convergence rate can be derived and the
% main results for the mixed FEM \refe{mixed-fem1}-\refe{mixed-fem2} are
% summarized in the following corollary.

\section{Estimates for problems with more regular boundary data}
\label{more-smooth-data}
\setcounter{equation}{0}

Classical mixed FEM theory covers the case $g \in H^{t}(\Gamma)$ with $t \ge 1/2$.
In this section, we investigate the model problem \refe{pde}-\refe{pde-bc}
with boundary data $g \in H^{t}(\Gamma)$ for $0 < t < 1/2$.
An improved convergence rate can be derived and the
main results for the mixed FEM \refe{mixed-fem1}-\refe{mixed-fem2} are
summarized in the following corollary.

\begin{corollary}
	\label{byproduct-theorem}
	
	Let $f \in L^2(\Omega)$, $g \in H^t(\Gamma)$ for $0 < t < 1/2$.
	The mixed FEM \refe{mixed-fem1}-\refe{mixed-fem2} admits a unique solution $u_h$
	satisfying the following estimate
	\begin{align}
		\| u_{h} - u \|_{L^2(\Omega)} \le
		C h^{t+s-\frac{1}{2}}\|g\|_{H^t(\Gamma)} + C h^{s} \|f\|_{L^2(\Omega)}\,,
		\label{byproduct-error-estimate-u}
	\end{align}
	where $u$ is the very weak solution.
\end{corollary}

\noindent{\bf Proof.}
The regularized boundary data $g^h$ satisfies
\begin{align}
	\|g - g^h\|_{H^{-s}(\Gamma)} \le C h^{t+s} \|g\|_{H^{t}(\Gamma)}
\end{align}
\textcolor{black}{From \cite[Remark 5.4]{Apel}}, we have an improved estimate for
the regularized solution $u^h$
\color{black}
\begin{align}
	\| u^h- u\|_{L^2(\Omega)}  \le  C h^{t+s-\frac{1}{2}} \|g\|_{H^{t}(\Gamma)}\,.
	\label{regularization-error-improved}
\end{align}
And the projection error satisfies
\begin{eqnarray}
	\|\P_h u^h - u^h\|_{L^2(\Omega)} & \le & C h \|u^h\|_{H^1(\Omega)}
	\nn
	\\
	& \le & C h (\|g^h\|_{H^{\frac{1}{2}}(\Gamma)} + \|f\|_{L^2(\Omega)})
	\nn
	\\
	& \le & C h^{t+\frac{1}{2}}\|g\|_{H^{t}(\Gamma)} + C h \|f\|_{L^2(\Omega)}\,.
	\label{error-uh-mid-byproduct}
\end{eqnarray}

Recall the proof of Lemma \ref{lemma-sigma-eror}, we have
\begin{align}
	\|\Pi_h\bs^h - \bs_h\|^2_{L^2(\Omega)}
	& \le C (\|\bs^h\|_{L^2(\Omega)} + h \|\D \bs^h\|_{L^2(\Omega)} )
	\|\Pi_h\bs^h -  \bs_h\|_{L^2(\Omega)}
	\nn                                                                     \\
	& \quad + Ch^{ t-\frac{1}{2} } \|g\|_{H^t(\Gamma)}
	\|(\Pi_h\bs^h-\bs_h) \cdot
	\bn\|_{H^{-\frac{1}{2}}(\Gamma)} 
	\nn                                                                     \\
	& \le C( h^{t-\frac{1}{2}}  \|g\|_{H^t(\Gamma)} + \|f\|_{L^2(\Omega)}) 
	\|\Pi_h\bs^h - \bs_h\|_{L^2(\Omega)}\,.
	\label{final-sigma-errorbyproduct}
\end{align}
and the quasi projection estimate
\begin{eqnarray}
	\|\Pi_h\bs^h - \bs^h\|_{L^2(\Omega)} & \le &  C \|g^h\|_{H^{\frac{1}{2}}(\Gamma)} + C\|f\|_{L^2(\Omega)}
	\nn
	\\
	& \le & C h^{t-\frac{1}{2}}\|g\|_{H^t(\Gamma)} + C\|f\|_{L^2(\Omega)}\,.
	\label{quasi-projection-sigma-error-pbyproduct}
\end{eqnarray}
An improved estimate for $\bs_h$ follows directly
\begin{align}
	\label{error-sigma-improved}
	\| \bs^h - \bs_h\|_{L^2(\Omega)} \le C( h^{t-\frac{1}{2}} \|g\|_{H^t(\Gamma)} + \|f\|_{L^2(\Omega)})\,.
\end{align}

Then, from \refe{error-u} the error  $\|\P_h u^h - u_h\|_{L^2(\Omega)}$ satisfies
\begin{align}
	\label{error-u-improved}
	\|\P_h u^h - u_h\|_{L^2(\Omega)}^2
	& = -(\D \bm{\omega}_h ,\P_h u^h - u_h)
	\nn                                                                                      \\
	\textrm{(by \refe{error1})} \quad
	& =  \langle g-g^h, \bm{\omega}_{h} \cdot \bn \rangle  -(\bs^h - \bs_h , \bm{\omega}_h)
	\nn                                                                                      \\
	& = \mathcal{J}_3 + \mathcal{J}_4\,,
\end{align}
where $\mathcal{J}_3$ can be bounded by
\begin{align}
	\label{error-J3}
	\left| \mathcal{J}_3 \right|
	& \le \|g-g^h\|_{L^2(\Gamma)} \|(\bm{\omega}_{h}- \Pi_h \nabla z) \cdot \bn \|_{L^2(\Gamma)}
	+ \|g-g^h\|_{L^2(\Gamma)} \|( \Pi_h \nabla z - \nabla z )
	\cdot \bn \|_{L^2(\Gamma)}
	\nn                                                                                                     \\[3.5pt]
	& \quad + \|g-g^h\|_{H^{-(s-\frac{1}{2})}(\Gamma)} \|\nabla z \cdot \bn \|_{H^{s-\frac{1}{2}}(\Gamma)}
	\nn                                                                                                     \\
	& \le C h^{t-\frac{1}{2}}\|(\bm{\omega}_{h}- \Pi_h \nabla z) \|_{L^2(\Omega)} \|g\|_{H^t(\Gamma)}
	+ C h^{t+s-\frac{1}{2}} \|\nabla z\|_{H^{s-\frac{1}{2}}(\Gamma)} \|g\|_{H^s(\Gamma)}
	\nn                                                                                                     \\[3.5pt]
	& \quad + C h^{t+s-\frac{1}{2}} \|g\|_{H^t(\Gamma)} \| z\|_{H^{s+1}(\Omega)}
	\nn                                                                                                     \\
	& \le C  h^{t+s-\frac{1}{2}}\|g\|_{H^{t}(\Gamma)}
	\|\P_h u^h -u_h\|_{L^2(\Omega)},
\end{align}
and $\mathcal{J}_4$ can be bounded by
\begin{align}
	\label{error-J4}
	\left| \mathcal{J}_4 \right|
	& \le \|\bs^h-\bs_h\|_{L^2(\Omega)} \|\bm{\omega}_h - \nabla z\|_{L^2(\Omega)}
	+ \|\D(\bs^h - \bs_h)\|_{L^2(\Omega)} \| z - \P_h z\|_{L^2(\Omega)}
	\nn                                                                                      \\[3.5pt]
	\textrm{(by \refe{error-sigma-improved} and \refe{standard-mixed-fems-error})} \quad
	& \le  C h^s  \left(h^{t-\frac{1}{2}}\|g\|_{H^t(\Gamma)} + \|f\|_{L^2(\Omega)} \right)
	\|\P_h u^h - u_h\|_{L^2(\Omega)}
	\nn
	\\
	& \quad  + Ch \|f\|_{L^2(\Omega)} \|z\|_{H^{1}(\Omega) }
	\nn                                                                                      \\
	& \le C\left(h^{t+s-\frac{1}{2}}\|g\|_{H^t(\Gamma)} + h^s\|f\|_{L^2(\Omega)} \right)
	\|\P_h u^h - u_h\|_{L^2(\Omega)} \,.
\end{align}
Taking the above two estimates into \refe{error-u-improved} yields an improved estimate
\begin{align}
	\|\P_h u^h - u_h\|_{L^2(\Omega)} \le
	C\left( h^{t+s-\frac{1}{2}}\|g\|_{H^t(\Gamma)} + h^s\|f\|_{L^2(\Omega)}  \right).
\end{align}
At last, combining the above estimates
with \refe{regularization-error-improved} and  \refe{error-uh-mid-byproduct},
Corollary \ref{byproduct-theorem} is proved.
\quad \endproof

\vspace{0.4cm}

%==============================
% Numerical Results
%==============================
\section{Numerical results}
\label{numericalresults}
\setcounter{equation}{0}

In this section, we provide several numerical examples to demonstrate the effectiveness of
the proposed mixed FEM \refe{mixed-fem1}-\refe{mixed-fem2}.
All computations are performed by the free software FEniCSx
\cite{BarattaEtal2023},
and the meshes are generated by Gmsh \cite{geuzaineGmsh3DFinite2009}.

\begin{example}
	\rm
	\label{example-convex}
	In the first example, we take a rectangular domain $\Omega=(-1,1) \times (0, 1)$.
	Then we consider the Poisson equation with Dirichlet boundary condition
	\begin{equation}
		\left\{\begin{aligned}
			-\Delta u & = 0 \qquad \text{in} \quad\Omega _\omega ,
			\\
			u         & = g  \qquad \text{on}  \quad\Gamma _\omega,
		\end{aligned}
		\right.
		\label{equ:numerical-problem}
	\end{equation}
	where, the exact solution in polar coordinates is defined by
	\begin{equation}
		u(r,\theta) = r^{-0.4999}\sin(-0.4999\theta).
		\label{artificial-solution}
	\end{equation}
	As $u$ is harmonic in $\Omega$, the source term $f=0$ belongs to $L^2(\Omega)$ and
	the boundary data $g$ can be simply defined as $g(r,\theta)=u(r,\theta) $ on $ \Gamma$.
	It is easy to verify that $g$ belongs to $L^2(\Gamma)$,  but $g$ is not in $H^{1/2}(\Gamma)$.
	
	\begin{figure}[htbp]
		\centering
		\includegraphics[width=0.24\linewidth]{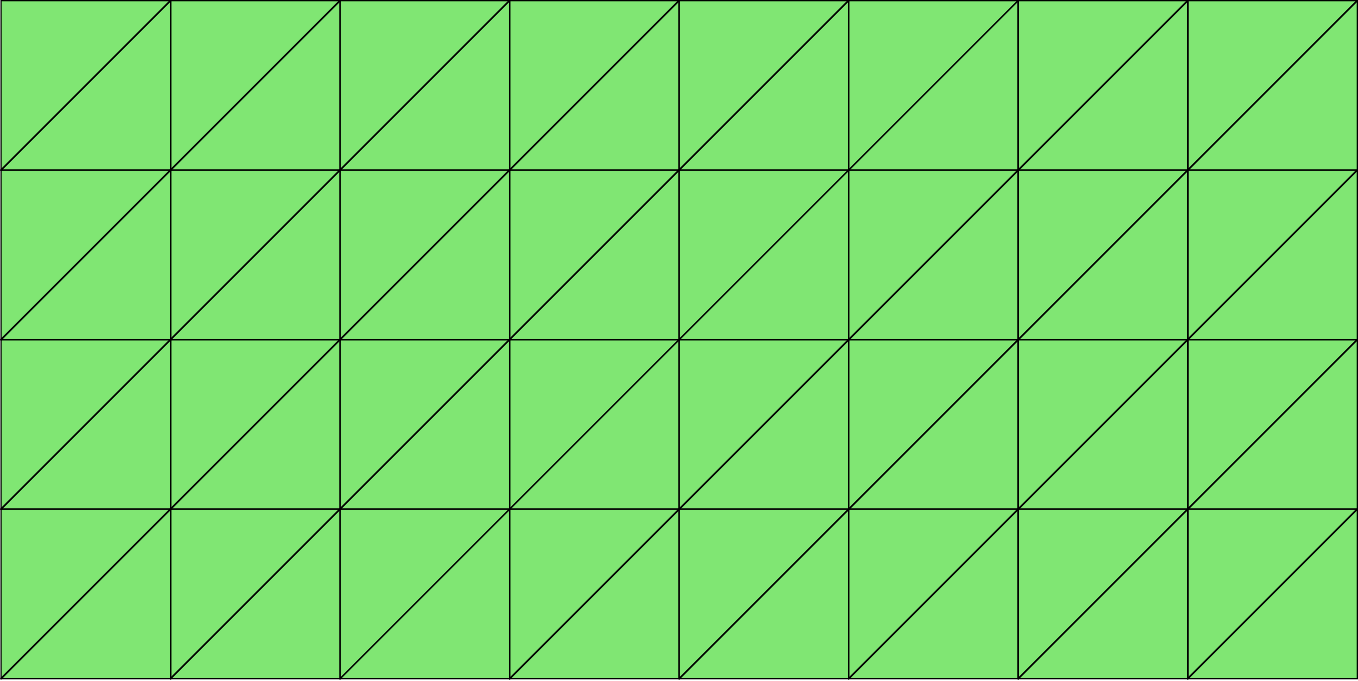}
		\caption{The uniform meshes with $h=\sqrt{2}/4$. (Example \ref{example-convex})}
		\label{fig:mesh-two-n4}
	\end{figure}
	
	We solve the above artificial problem by the proposed mixed FEM \refe{mixed-fem1}-\refe{mixed-fem2}
	on uniform triangular meshes, see Figure \ref{fig:mesh-two-n4} for an illustration.
	The plots of $u_h$ with $h=\sqrt{2}/32$ are shown in Figure \ref{fig:mesh-two-n4}.
	For comparison, we use the standard linear Lagrange FE with $L^2(\Gamma)$ projection on the boundary
	to solve this artificial problem.
	The numerical $u_h$ computed by the linear FE  $P_1$ on the same mesh with
	$h=\sqrt{2}/32$ is also shown in Figure \ref{fig:solution-square-n32}.
	We can observe numerical oscillation near the singular boundary points.
	
	\begin{figure}[htbp]
		\centering
		\includegraphics[width=0.49\linewidth]{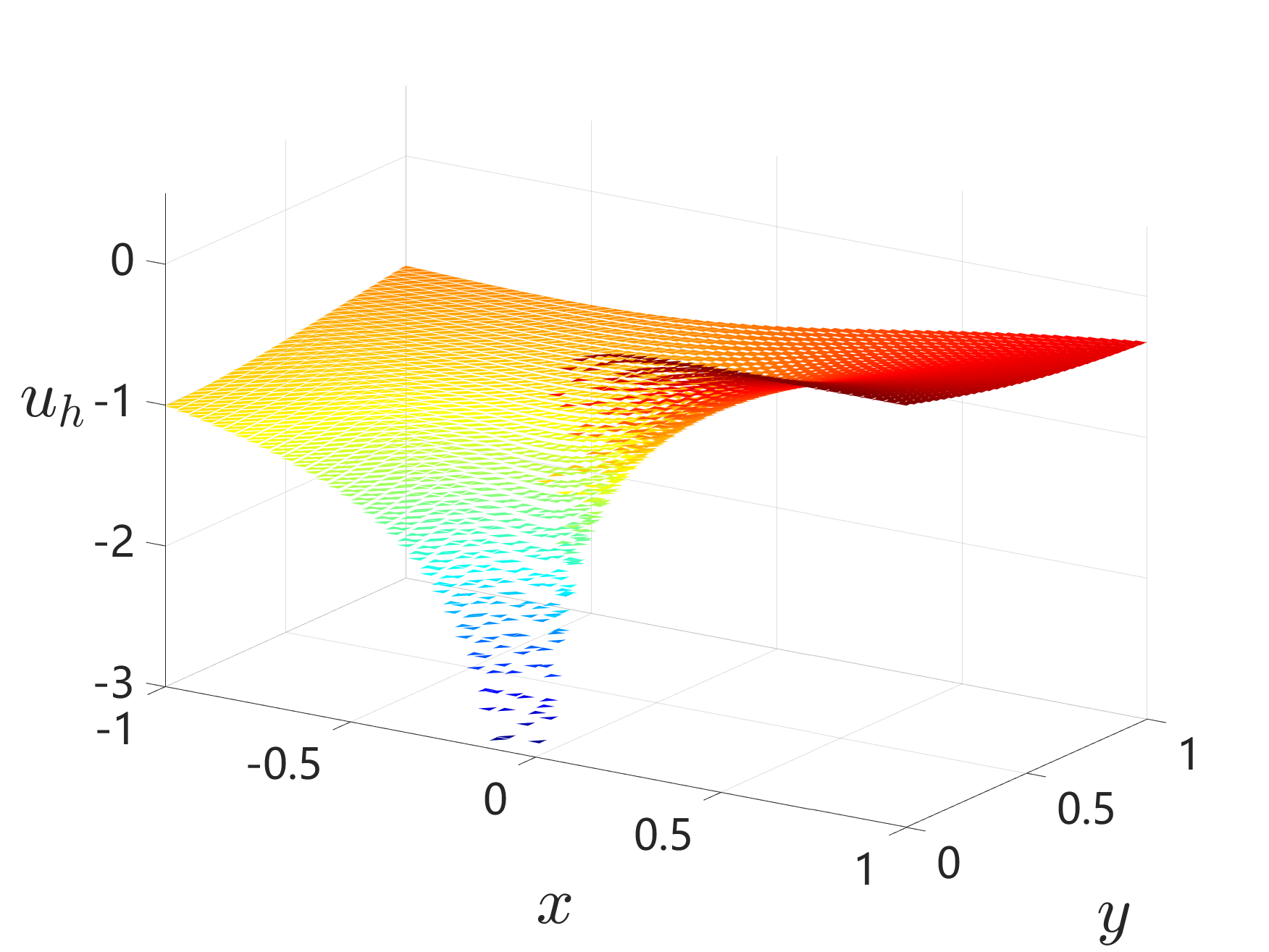}
		\includegraphics[width=0.49\linewidth]{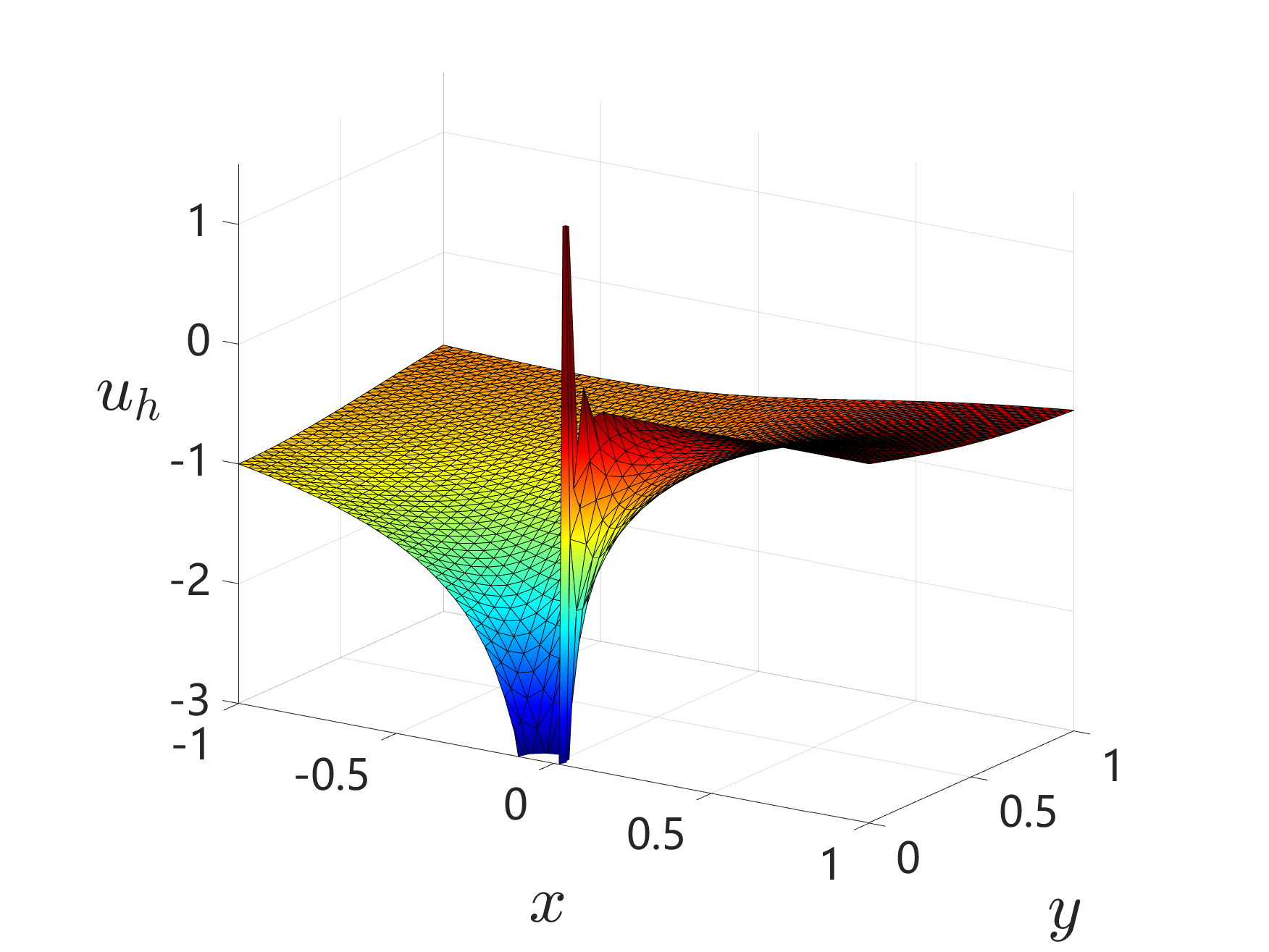}
		\caption{Numerical $u_h$ computed by mixed FEM $\RT_{\! \! 0} \times \DG_{\! 0}$(Left);
			Numerical $u_h$ computed by linear FEM $P_1$ with $L^2(\Gamma)$ projection(Right).
			(Example \ref{example-convex})}
		\label{fig:solution-square-n32}
	\end{figure}
	The $L^2$-norm errors $\|u-u_h\|_{L^2(\Omega)}$ on gradually refined meshes
	are presented in Table \ref{tab:convex},
	which clearly shows the $O(h^{1/2})$ convergence.
	Moreover, we also provide the errors $\|\bs_h-\bs^h\|_{L^2(\Omega)}$.
	Since $\bs^h \in \bH(\D)$ is unknown,
	based on the regularized boundary data $g^h$
	we compute an approximation of $\bs^h$ by a linear FEM on
	a fine mesh with mesh size $\sqrt{2}/{512}$.
	The errors in Table \ref{tab:convex} indicate that
	an $O(h^{-1/2})$ convergence for  $\|\bs_h-\bs^h\|_{L^2(\Omega)}$, which implies that
	the estimate for $\bs_h$ in \refe{error-sigma} is sharp.
	
	\begin{table}[htbp]
		\centering
		\caption{Errors and convergence rates for the rectangular domain. (Example \ref{example-convex})}
		\label{tab:convex}
		\begin{tabular*}{\hsize}{@{}@{\extracolsep{\fill}}cccccc@{}}
			\toprule
			$h$   & $\|u-u_h\|_{L^2(\Omega)}$ & Rate     & $\|\bs_h-\bs^h\|_{L^2(\Omega)}$ & Rate      \\
			\midrule
			$\sqrt{2}/2  $ & 0.335280    & ---      & 2.119086          & ---       \\
			$\sqrt{2}/4  $ & 0.244516    & 0.455435 & 2.994347          & -0.498799 \\
			$\sqrt{2}/8  $ & 0.175349    & 0.479701 & 4.236726          & -0.500709 \\
			$\sqrt{2}/16 $ & 0.124972    & 0.488626 & 5.997160          & -0.501330 \\
			$\sqrt{2}/32 $ & 0.088831    & 0.492463 & 8.508301          & -0.504591 \\
			$\sqrt{2}/64 $ & 0.063064    & 0.494245 & 12.160640         & -0.515276 \\
			$\sqrt{2}/128$ & 0.044745    & 0.495109 & 17.766272         & -0.546922 \\
			\bottomrule
		\end{tabular*}
	\end{table}
	
\end{example}

\begin{example}
	\label{example-concave}
	\rm
	
	\begin{figure}[htbp]
		\centering
		\includegraphics[width=0.2\linewidth]{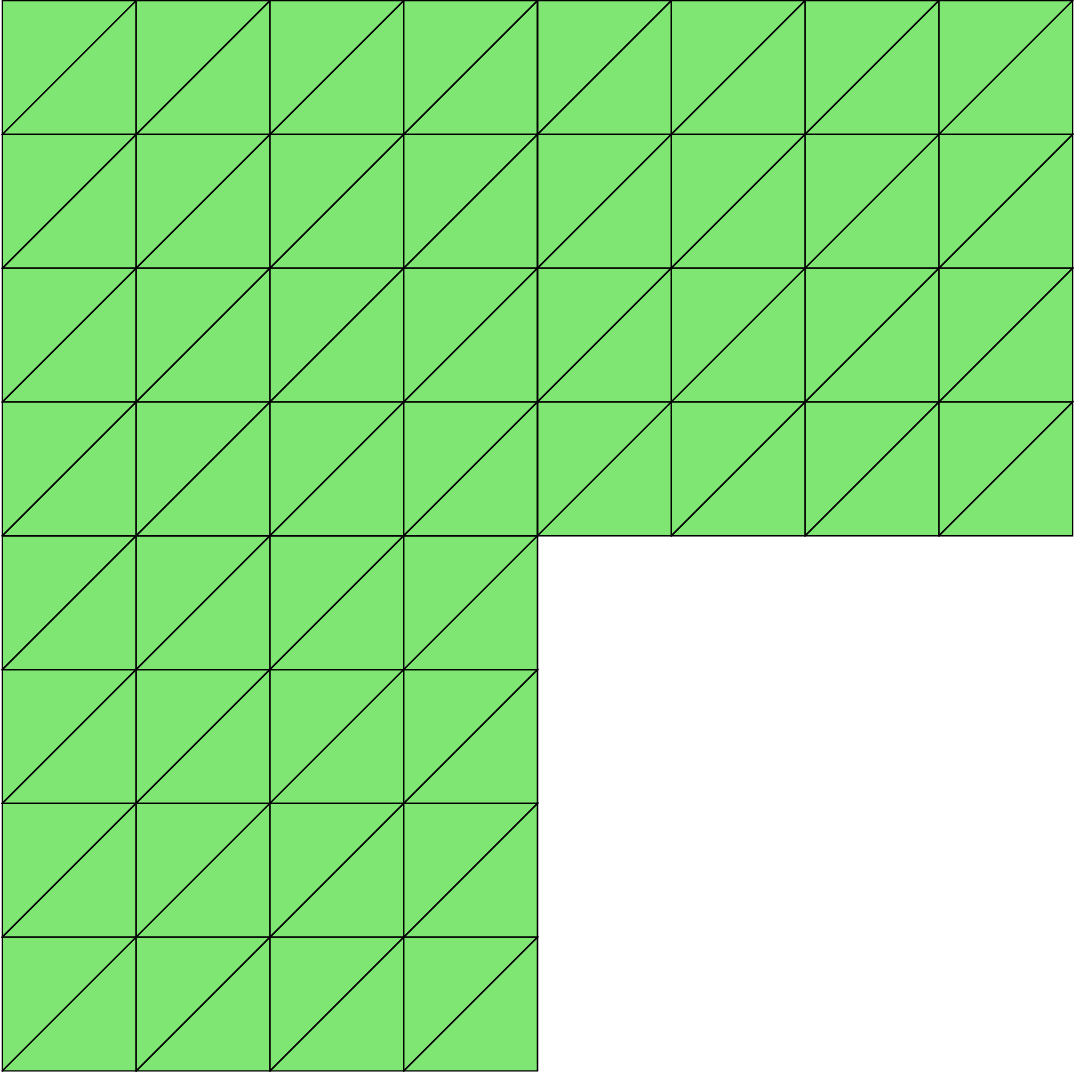}
		\caption{The uniform meshes with $h=\sqrt{2}/4$.(Example \ref{example-concave})}
		\label{fig:mesh-lshape-n4}
	\end{figure}
	
	In the second example, we solve the problem \refe{equ:numerical-problem} in
	a nonconvex L-shape domain $\Omega=(-1,1)^2-[0,1) \times (-1,0]$.
	We take the same exact solution $u$ in \refe{artificial-solution} and the boundary data $g=u|_{\Gamma}$.
	A uniform mesh is used in our tests, see Figure \ref{fig:mesh-lshape-n4} for illustration.
	The numerical results $u_h$ with $h=\sqrt{2}/{32}$ are plot in Figure \ref{fig:solution-lshape-n32}.
	For comparison, we also show the plot of $u_h$ computed on the same mesh
	by conventional linear FEM with $L^2(\Gamma)$ projection
	in Figure \ref{fig:solution-lshape-n32}.
	Again, we observe the numerical oscillation near the singular point,
	which agrees with previous numerical results in \cite[Figure 2]{Apel}.

	\begin{figure}[hbtp]
		\centering
		\includegraphics[width=0.49\linewidth]{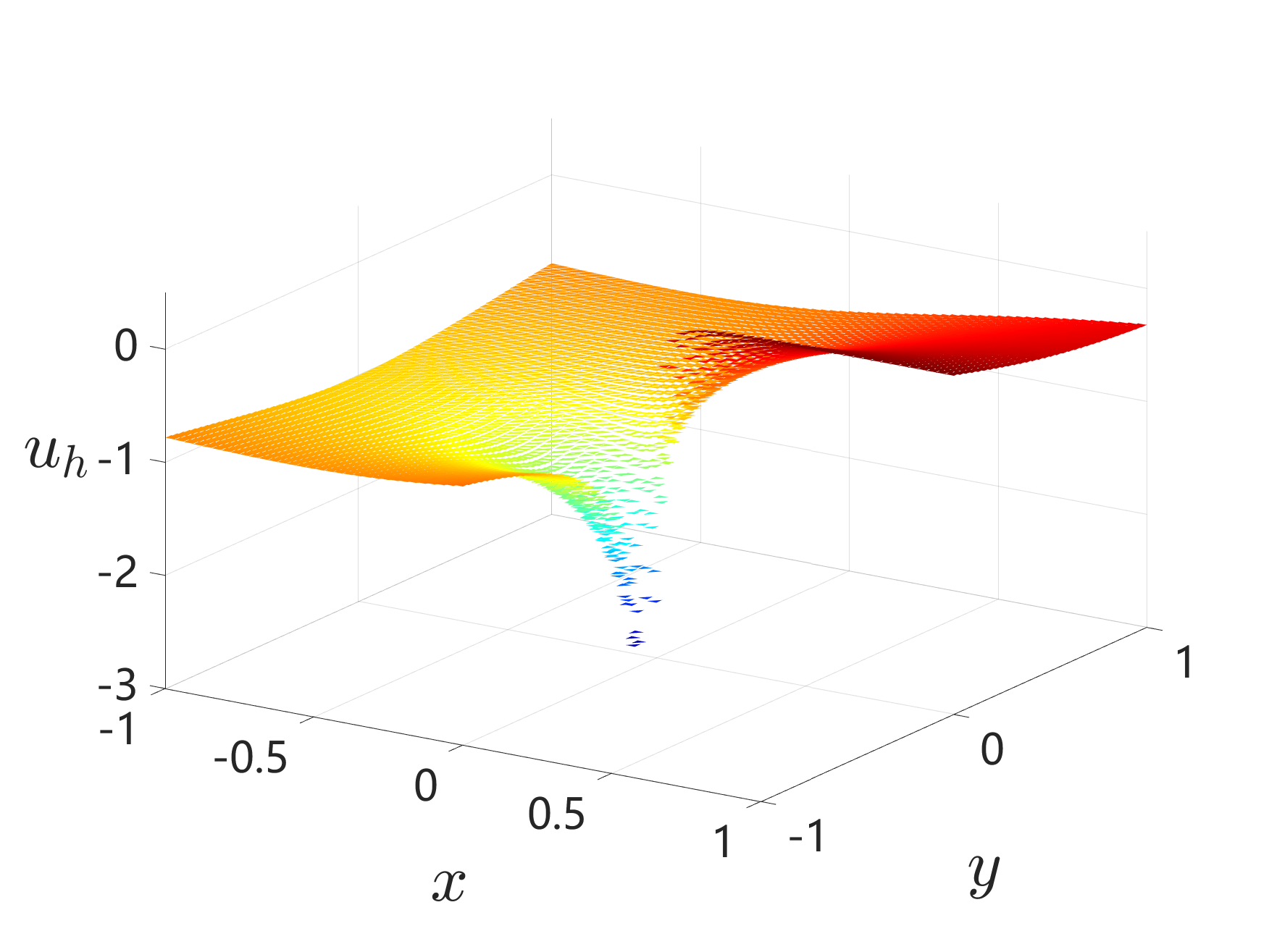}
		\includegraphics[width=0.49\linewidth]{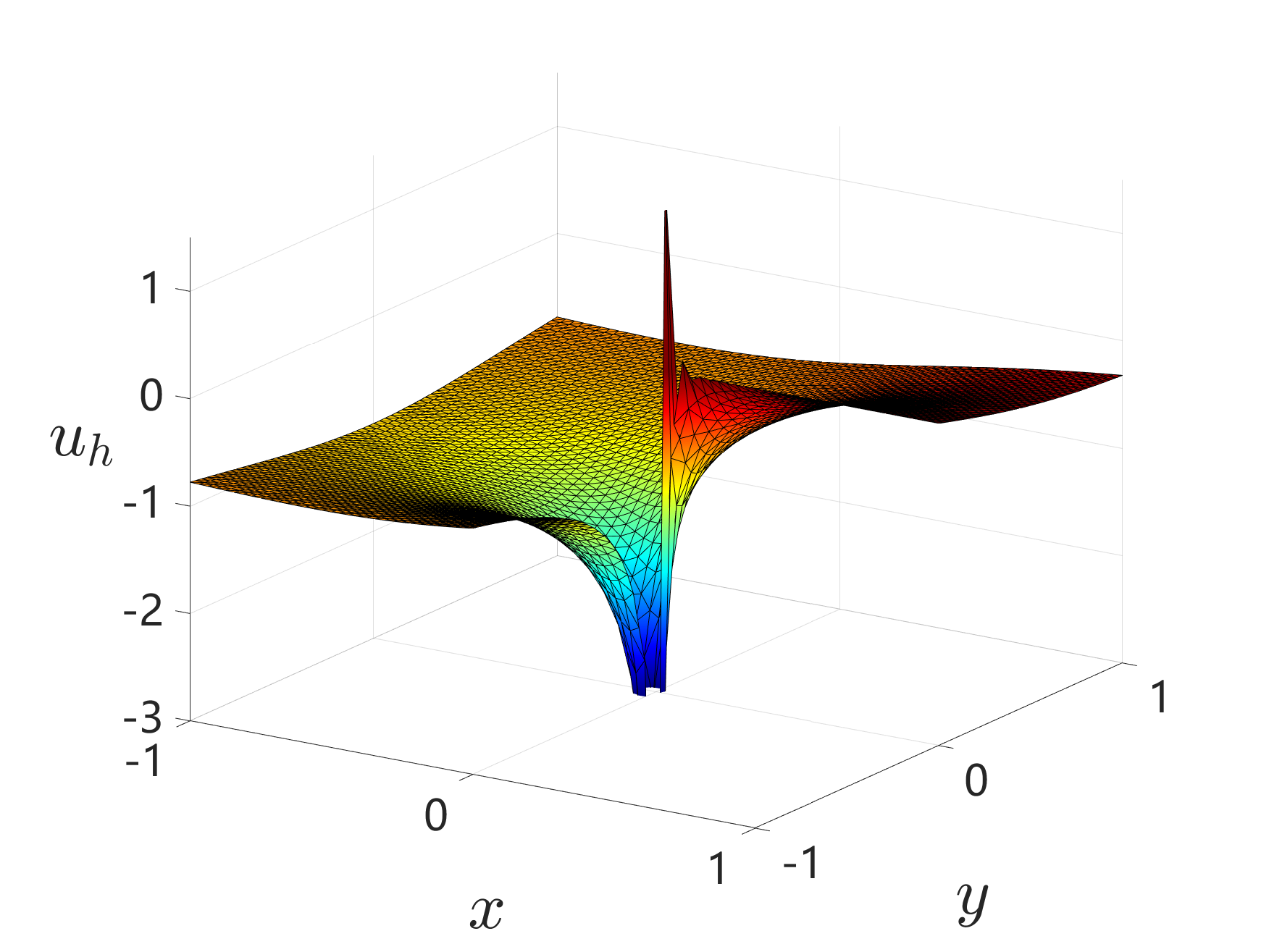}
		\caption{Numerical $u_h$ computed by mixed FEM $\RT_{\! \! 0} \times \DG_{\! 0}$(Left);
			Numerical $u_h$ computed by linear FEM $P_1$ with $L^2(\Gamma)$ projection(Right).
			(Example \ref{example-concave})}
		\label{fig:solution-lshape-n32}
	\end{figure}

	We show the errors $\|u-u_h\|_{L^2(\Omega)}$ and $\|\bs_h-\bs^h\|_{L^2(\Omega)}$
	on gradually refined meshes in Table \ref{tab:non-convex}.
	As $\Theta={3\pi}/{2}$ for the L-shape domain,
	the convergence rate for $\|u-u_h\|_{L^2(\Omega)}$ is nearly $O(h^{1/6})$
	according to estimate \refe{main-error-estimate-u} in Theorem \ref{main-theorem}.
	One can easily observe that the error results agree with our theoretical results well.
	
	\begin{table}[htbp]
		\centering
		\caption{Errors of $(\bs_h,u_h)$ in the L-shape domain.(Example \ref{example-concave})}
		\label{tab:non-convex}
		\begin{tabular*}{\hsize}{@{}@{\extracolsep{\fill}}cccccc@{}}
			\toprule
			h     & $\|u-u_h\|_{L^2(\Omega)}$ & Rate     & $\|\bs_h-\bs^h\|_{L^2(\Omega)}$ & Rate      \\
			\midrule
			$\sqrt{2}/2  $ & 0.681983    & ---      &   1.540822       & ---       \\
			$\sqrt{2}/4  $ & 0.598987    & 0.187213 &   2.214954       &  -0.523577 \\
			$\sqrt{2}/8  $ & 0.525100    & 0.189931 &   3.161334       &  -0.513256 \\
			$\sqrt{2}/16 $ & 0.461639    & 0.185828 &   4.497357       &  -0.508544 \\
			$\sqrt{2}/32 $ & 0.407324    & 0.180590 &   6.410924       &  -0.511455 \\
			$\sqrt{2}/64 $ & 0.360495    & 0.176196 &   9.227381       &  -0.525389 \\
			$\sqrt{2}/128$ & 0.319760    & 0.172990 &  13.617184       &  -0.561435 \\
			\bottomrule
		\end{tabular*}
	\end{table}

\end{example}

\begin{example}
	\rm
	\label{example-improved}
	In the third example, we consider the Poisson equation with
	boundary data $g \in H^{s}(\Gamma)$ with $0 < s < 1/2$.
	Here, we take the exact solution
	$u=r^{-\frac{1}{3}}\sin(-\frac{1}{3} \theta)$ in polar coordinates,
	where the boundary data $g \in H^{t}(\Gamma)$ for any $t <1/6$.
	
	We test the performance of the mixed FEM \refe{mixed-fem1}-\refe{mixed-fem2} for
	both rectangular and L-shape domains.
	For the rectangular domain which is convex,
	The error estimates in Corollary \ref{byproduct-theorem} indicate that
	$ \| u_{h} - u \|_{L^2(\Omega)} $ is around $O(h^{2/3})$.
	The numerical errors for the rectangular domain in Table \ref{tab:nonconvex:Hs}
	agree with our theoretical results.
	For the L-shape domain,   Corollary \ref{byproduct-theorem} implies that
	$ \| u_{h} - u \|_{L^2(\Omega)} $ is around $O(h^{1/3})$.
	The numerical results in Table \ref{tab:convex:Hs} clearly show that
	our estimate is sharp.

	\begin{table}[htbp]
		\centering
		\caption{Errors in the rectangular domain with $g \in H^{1/6-\epsilon}(\Gamma)$. (Example \ref{example-improved})}
		\label{tab:convex:Hs}
		\begin{tabular*}{\hsize}{@{}@{\extracolsep{\fill}}cccccc@{}}
			\toprule
			h     & $\|u-u_h\|_{L^2(\Omega)}$ & Rate     & $\|\bs_h-\bs^h\|_{L^2(\Omega)}$ & Rate      \\
			\midrule
			$\sqrt{2}/2  $ & 0.151589   & ---      & 1.066496    & ---       \\
			$\sqrt{2}/4  $ & 0.100904   & 0.587177 & 1.343957    & -0.333608 \\
			$\sqrt{2}/8  $ & 0.065459   & 0.624334 & 1.694723    & -0.334563 \\
			$\sqrt{2}/16 $ & 0.041955   & 0.641744 & 2.137200    & -0.334673 \\
			$\sqrt{2}/32 $ & 0.026712   & 0.651351 & 2.700412    & -0.337457 \\
			$\sqrt{2}/64 $ & 0.016941   & 0.657005 & 3.434433    & -0.346893 \\
			$\sqrt{2}/128$ & 0.010718   & 0.660436 & 4.454866    & -0.375310 \\
			\bottomrule
		\end{tabular*}
	\end{table}
	
	\begin{table}[htbp]
		\centering
		\caption{Errors in the L-shape domain with $g \in H^{1/6-\epsilon}(\Gamma)$.(Example \ref{example-improved})}
		\label{tab:nonconvex:Hs}
		\begin{tabular*}{\hsize}{@{}@{\extracolsep{\fill}}cccccc@{}}
			\toprule
			h     & $\|u-u_h\|_{L^2(\Omega)}$ & Rate     & $\|\bs_h-\bs^h\|_{L^2(\Omega)}$ & Rate      \\
			\midrule
			$\sqrt{2}/2  $ & 0.284134   & ---      & 1.267748      & ---       \\
			$\sqrt{2}/4  $ & 0.212401   & 0.419782 & 1.604859      & -0.340179 \\
			$\sqrt{2}/8  $ & 0.159163   & 0.416283 & 2.029146      & -0.338426 \\
			$\sqrt{2}/16 $ & 0.120545   & 0.400940 & 2.564409      & -0.337754 \\
			$\sqrt{2}/32 $ & 0.092398   & 0.383641 & 3.249484      & -0.341584 \\
			$\sqrt{2}/64 $ & 0.071562   & 0.368668 & 4.153913      & -0.354260 \\
			$\sqrt{2}/128$ & 0.055866   & 0.357226 & 5.434413      & -0.387653 \\
			\bottomrule
		\end{tabular*}
	\end{table}
	
\end{example}

\vspace{0.4cm}

%==============================
%  Section conclusions
%==============================
\section{Conclusions}
\label{sec-conclusion}

In this paper,
we have extended the applicability of the Raviart--Thomas mixed method by rigorously proving that
it is suitable for solving elliptic problems with rough Dirichlet boundary data.
To the best of our knowledge, no analysis has been established for
the Raviart--Thomas mixed FEM for solving problems with boundary data in $L^2(\Gamma)$ only.
More important is that the Raviart--Thomas mixed FEM does not need to modify the boundary data,
although our proof is based on a regularized approach.
Numerical experiments presented in this work demonstrate the efficiency of the method and
confirm our theoretical analysis.

In this work, we confine our study to rough boundary data problems.
It is assumed that the source $f \in L^2(\Omega)$.
However, our results can be applied to the case $f \in H^{-1}(\Omega)$,
see \cite{GjerdeKN2021},  where Gjerde et al. use mixed FEM
to solve Poisson's problems with line sources.
As the exact solution $u \notin H^1(\Omega)$,
we only consider the lowest order mixed FEM $\RT_{\! \! 0} \times \DG_{\! 0}$.
Moreover, adaptive meshes might improve the performance of the numerical methods,
in particular for the nonconvex domain.

\color{black}
\section*{Declarations}

\noindent {\bf The Conflict of Interest Statement:}
No conflict of interest exists.
\vskip 0.1in

\noindent {\bf Availability of data and material:}
The code to reproduce the numerical results presented in this paper is available at
{\color{black}\verb|https://github.com/bombeuler/Mixed-FEM-Codes|}.
\vskip 0.1in

%==============================
%  The references
%==============================

\end{document}